\documentclass[reqno,b5paper]{amsart}
\usepackage{amsmath}
\usepackage{amssymb}
\usepackage{amsthm}
\usepackage{enumerate}
\usepackage[mathscr]{eucal}
\usepackage{epsfig}

\newtheorem{theo}{Theorem}[section]

\newtheorem{Cor}{Corollary}[section]

\newfont{\sss}{msam10 scaled\magstep0}
\newfont{\bbb}{msbm10 scaled\magstep1}

\def\R{\mbox{\bbb R}}
\def\Q{\mbox{\bbb Q}}
\def\N{\mbox{\bbb N}}

\begin{document}

\title[Fixed Points with an Iterate at a Point]{Fixed Points for Mappings with a Contractive Iterate at Each Point}
\author[SAMET KARAIBRYAMOV AND BOYAN ZLATANOV]{SAMET KARAIBRYAMOV* AND BOYAN ZLATANOV**}

\newcommand{\acr}{\newline\indent}
\address{\llap{*\,}Department of Mathematics and Informatics \acr
                              of Mathematics and Informatics \acr
                              Plovdiv University \acr
                              24 ``Tzar Assen'' str. \acr
                              Plovdiv, 4000 \acr
                              BULGARIA}
\email{skaraibryamov@gmail.com}

\address{\llap{**\,}Department of Mathematics and Informatics \acr
                              of Mathematics and Informatics \acr
                              Plovdiv University \acr
                              24 ``Tzar Assen'' str. \acr
                              Plovdiv, 4000 \acr
                              BULGARIA}
\email{bzlatanov@gmail.com}
\thanks{Research is partially supported by Plovdiv University ``Paisii Hilendarski'', NPD, Project NI11--FMI-004.}
\begin{abstract}
We generalize the results of Sehgal and Guseman for
mappings on a complete metric space with a contractive iterate condition at each point.
\end{abstract}

\maketitle

\section{Introduction}

Let $\bigl(X,\rho \bigr)$ be a metric space and $T:X\to X$ be a map. If there is $\alpha <1$, such that
$\rho (Tx,Ty)\leq \alpha\rho (x,y)$
holds for every $x,y\in X$, then $T$ is called a contraction. A well known theorem of Banach \cite{B}
states that if $X$ is a complete and $T$ is a contraction, then $T$ has a unique fixed point $z$ and for any $x\in X$
the sequence of the successive approximations $\bigl\{T^nx\bigr\}_{n=1}^\infty$ converges to $z$, where $T^nx=T\bigl(T^{n-1}x\bigr)$.
Sehgal has generalized this result in \cite{S}, by considering maps
that a contractive at a point. The result of Sehgal has been generalized by Guseman \cite{G}. In
these two papers the authors have imposed weaker contraction conditions on the map, but
again the fixed point can always be found by using the Picard iteration, beginning with some
initial choice $x\in X$. Some other articles that consider mappings that are with a contrective
iterate at a point are \cite{C2, C, GL, GL2, M}.

\section{Main result}
\begin{theo}\label{main}
Let $\bigl(X,\rho \bigr)$ be a complete metric space and $T:X\to X$ be a map with the properties:
\renewcommand{\theenumi}{\textnormal{(\alph{enumi})}}
\renewcommand{\labelenumi}{\theenumi}
\begin{enumerate}
  \item there exist subsets $U_n\subseteq X$, $n\in\N$, such that $T:U_n\to U_n$, and $\bigcap_{n=1}^\infty U_n\not=\emptyset$;
  \item there exists $\alpha\in (0,1)$, such that for every $x\in X$ there are $n_1, n_0\in\N$, such that for every $u\in \bigcup_{n=n_0}^\infty U_n$ there holds the inequality
$$
\rho (T^{n_1}u,T^{n_1}x)\leq \alpha\rho (u,x).
$$
\end{enumerate}
Then there exists $z\in X$ with the properties:
\renewcommand{\theenumi}{\textnormal{(\roman{enumi})}}
\renewcommand{\labelenumi}{\theenumi}
\begin{enumerate}
  \item There exists $n_z\in\N$, such that $T^{n_z}z=z$;\label{vi}
  \item For every $y\in X$ there is a sequence of naturals
$\{q_i\}_{i=1}^\infty$, such that $\displaystyle\lim_{i\to\infty}T^{q_i}y=z$;\label{vii}
  \item For every $x\in \bigcap_{n=1}^\infty U_n$ there holds $\lim_{n\to\infty}T^nx=z$.\label{viii}
\end{enumerate}
\end{theo}
\begin{proof}
Let $x\in\bigcap_{n=1}^\infty U_n$ be arbitrary chosen.
By $T:U_n\to U_n$ it follows that for every $s\in\N$ we have $T^sx\in\bigcap_{n=1}^\infty U_n$.
Following \cite{G} and \cite{S} let define the function $r(x)=\sup\{\rho (T^nx,x):n\in\N\}$ for any $x\in X$. First we will prove that for any $x\in\bigcap_{n=1}^\infty U_n$ there holds the inequality $r(x)<+\infty$.

Let choose an arbitrary $x\in\bigcap_{n=1}^\infty U_n$, then there exist $n_1(x), n_0(x)\in\N$, such that the inequality $\rho (T^{n_1(x)}u,T^{n_1(x)}x)\leq \alpha\rho (u,x)$
holds for any $u\in\bigcup_{n=n_0(x)}^\infty U_n$. Therefore for any $s\in\N$ there holds
the inequality
$$
\rho (T^{n_1(x)}(T^sx),T^{n_1(x)}x)\leq \alpha\rho (T^sx,x).
$$
Put $l(x)=\max\{\rho (T^sx,x),s=1,2,\dots ,n_1\}$.
For any $s\in\N$ there is $k\in\N$, such that $kn_1(x)<s\leq (k+1)n_1(x)$, then we can write the chain of inequalities:
$$
\begin{array}{lll}
\rho (T^sx,x)&\leq&\rho (T^{n_1(x)}(T^{s-n_1(x)}x),T^{n_1(x)}x)+\rho (T^{n_1(x)}x,x)\\[12pt]
&\leq&\alpha\rho (T^{s-n_1(x)}x,x)+\rho (T^{n_1(x)}x,x)\leq\alpha\rho (T^{s-n_1(x)}x,x)+l(x)\\[12pt]
&\leq&\alpha\rho (T^{n_1(x)}(T^{s-2n_1(x)}x),T^{n_1(x)}x)+\alpha\rho (T^{n_1(x)}x,x)+l(x)\\[12pt]
&\leq&\alpha^2\rho (T^{s-2n_1(x)}x,x)+(1+\alpha)l(x)\\[12pt]
&\leq&\alpha^2\rho (T^{n_1(x)}(T^{s-3n_1(x)}x),T^{n_1(x)}x)+\alpha^2\rho (T^{n_1(x)}x,x)\\[12pt]
&&+(1+\alpha)l(x)\\[12pt]
&\leq&\alpha^3\rho (T^{s-3n_1(x)}x,x)+(1+\alpha+\alpha^2)l(x)\\[12pt]
&\leq&\dotfill\\[12pt]
&\leq&\alpha^k\rho (T^{s-kn_1(x)}x,x)+(1+\alpha+\alpha^2+\cdots +\alpha^{k-1})l(x)\\[12pt]
&\leq& l(x)\displaystyle\sum_{j=0}^k \alpha^j \leq\displaystyle\frac{l(x)}{1-\alpha}.
\end{array}
$$
Thus for any $x\in\bigcap_{n=1}^\infty U_n$ there holds the inequality $r(x)\leq\displaystyle\frac{l(x)}{1-\alpha }<+\infty$.

Now we will construct inductively a convergent sequence $\{x_i\}_{i=1}^\infty$.
Let choose an arbitrary $x\in\bigcap_{n=1}^\infty U_n$.
Put $m_0=n_1(x)$, $x_1=T^{m_0}x$, $m_1=n_1(x_1)$,
$x_2=T^{m_1}x_1$. If we have defined $m_{i-1}\in\N$ and $x_{i}\in\bigcap_{n=1}^\infty U_n$, then put $m_i=n_1(x_i)$ and $x_{i+1}=T^{m_i}x_i$.

We will show that $\{x_i\}_{i=1}^\infty$ is a Cauchy sequence. Indeed by the inequalities
$$
\begin{array}{lll}
\rho (x_{i+1},x_i)&=&\rho (T^{m_i}x_i,x_i)=\rho (T^{m_{i-1}}(T^{m_i}x_{i-1}),T^{m_{i-1}}x_{i-1})\\[12pt]
&\leq&\rho (T^{m_i}x_{i-1},x_{i-1})=\alpha\rho (T^{m_{i-2}}(T^{m_i}x_{i-2}),T^{m_{i-2}}x_{i-2})\\[12pt]
&\leq& \alpha^2\rho (T^{m_i}x_{i-2},x_{i-2})\leq\dots\leq \alpha^i\rho (T^{m_i}x,x)\leq \alpha^ir(x)\\[12pt]
\end{array}
$$
we get
$$
\begin{array}{lll}
\rho (x_{i+p},x_i)&=&\rho (x_{i+p},x_{i+p-1})+\rho (x_{i+p-1},x_{i+p-2})+\cdots\\[12pt]
&&\cdots+\rho (x_{i+2},x_{i+1})+\rho (x_{i+1},x_{i})\\[12pt]
&\leq& r(x)\sum_{s=i}^{i+p-1}\alpha^i\leq r(x)\frac{\alpha^i}{1-\alpha},
\end{array}
$$
which ensures that $\{x_i\}_{i=1}^\infty$ is a Cauchy sequence and taking into account the completeness of $X$ it follows that
there is $z\in X$, such that $\lim_{i\to\infty}x_i=z$.

Let for the rest of the proof $\{x_i\}_{i=1}^\infty$ be the sequence, that was constructed above and $z$
 be the limit of the sequence $\{x_i\}_{i=1}^\infty$.

We will show that $\displaystyle\lim_{s\to\infty}\rho (T^{s}x,z)=0$.
Choose $\varepsilon >0$. There is $i_0\in\N$, such that the inequalities $\rho (z,x_i)<\varepsilon /2$ and $\alpha^{i}r(x)<\varepsilon /2$ hold for every $i\geq i_0$.
Put $s_0=\sum_{i=0}^{i_0}m_i$. Then for every $s>s_0$
we can write the chain of inequalities
$$
\begin{array}{lll}
\rho (z,T^sx)&\leq&\rho (z,T^{m_{i_0}}x_{i_0})+\rho (T^{m_{i_0}}x_{i_0},T^sx)\\[12pt]
&\leq&\rho (z,x_{i_0+1})+\alpha\rho (x_{i_0},T^{s-m_{i_0}}x)\\[12pt]
&=&\displaystyle\frac{\varepsilon}{2}+\alpha\rho (T^{m_{i_0-1}}x_{i_0-1},T^{s-m_{i_0}}x)\\[12pt]
&\leq&\displaystyle\frac{\varepsilon}{2}+\alpha^2\rho (x_{i_0-1},T^{s-(m_{i_0}+m_{i_0-1})}x)\\[12pt]
&\leq& \dots\leq\displaystyle\frac{\varepsilon}{2}+\alpha^{i_0+1}\rho (x,T^{s-s_0}x)\leq\displaystyle\frac{\varepsilon}{2}+\alpha^{i_0+1}r(x)<\varepsilon .\\[12pt]
\end{array}
$$
By the arbitrary choice of $\varepsilon >0$ it follows that $\lim_{s\to\infty}\rho (T^{s}x,z)=0$.

There are $n_1(z), n_0(z)\in\N$, such that $\rho (T^{n_1(z)}u,T^{n_1(z)}z)\leq \alpha\rho (u,z)$, holds for any $u\in\bigcup_{n=n_0(z)}^\infty U_n$. By the construction of the sequence $\{x_i\}_{i=1}^\infty$
it follows that for any $s,i\in\N$ there holds $T^sx_i\in\bigcap_{n=1}^\infty U_n$ and thus the inequality
\begin{equation}\label{eq-2}
\rho (T^{n_1(z)}x_i,T^{n_1(z)}z)<\alpha\rho (x_i,z)
\end{equation}
holds for any $i\in\N$. Let $\varepsilon >0$ be arbitrary chosen, then there exists $j_0\in\N$, such that for any
$i\geq j_0$ there hold $\rho (x_i,z)<\varepsilon /3$ and $\alpha^ir(x)<\varepsilon /3$. Consequently by (\ref{eq-2}) we get
$\rho (T^{n_1(z)}x_i,T^{n_1(z)}z)<\varepsilon /3$ and the chain of inequalities
$$
\begin{array}{lll}
\rho (T^{n_1(z)}x_i,x_i)&=&\rho (T^{m_{i-1}}(T^{n_1(z)}x_{i-1}),T^{m_{i-1}}x_{i-1})\\[12pt]
&\leq& \alpha\rho (T^{n_1(z)}x_{i-1},x_{i-1})\\[12pt]
&=&\alpha\rho (T^{m_{i-2}}(T^{n_1(z)}x_{i-2}),T^{m_{i-2}}x_{i-2})\\[12pt]
&\leq& \alpha^2\rho (T^{n_1(z)}x_{i-2},x_{i-2})\\[12pt]
&\leq&\dots\leq \alpha^i\rho (T^{n_1(z)}x,x)\leq \alpha^ir(x)<\varepsilon /3.
\end{array}
$$

We will proceed with proving of \ref{vi}--\ref{viii}.

\ref{vi} Indeed we have that for any $\varepsilon >0$ there exists $j_0\in\N$, such that the inequality
$$
\rho (T^{n_1(z)}z,z)\leq \rho (T^{n_1(z)}z,T^{n_1(z)}x_i)+\rho (T^{n_1(z)}x_i,x_i)+\rho (x_i,z)<\varepsilon
$$
holds for every $i\geq j_0$ and thus by the arbitrary choice of $\varepsilon >0$ it follows that $T^{n_1(z)}z=z$.

\ref{vii} Let $y\in X$ be arbitrary.
We will construct inductively a sequence $\{y_i\}_{i=1}^\infty$. Put $p_0=n_1(y)$, $y_1=T^{p_0}y$, $p_1=n_1(y_1)$,
$y_2=T^{p_1}y_1$. If we have defined $p_{i-1}\in\N$ and $y_{i}$ then put $p_i=n_1(y_i)$,
$y_{i+1}=T^{p_i}y_i$ and $q_i=\sum_{j=0}^ip_j$.
By the construction of the sequence $\{x_i\}_{i=1}^\infty$ it follows that for any $s,i\in\N$ there holds
$T^sx_i\in\bigcap_{n=1}^\infty U_n$, thus
$$
\begin{array}{lll}
\rho (T^{q_i}x,T^{p_i}y_i)&\leq& \alpha\rho (T^{q_{i-1}}x,y_i)=\alpha\rho (T^{q_{i-1}}x,T^{p_{i-1}}y_{i-1})\\[12pt]
&\leq& \alpha^2\rho (T^{q_{i-2}}x,y_{i-1})\leq\dots\leq \alpha^{i+1}\rho (x,y).
\end{array}
$$
Therefore $\lim_{i\to\infty}\rho (T^{q_i}x,T^{p_i}y_i)=0$. Since $\lim_{s\to\infty}\rho (T^{s}x,z)=0$ and the triangle inequality
$$
\rho (T^{p_i}y_i,z)\leq\rho (T^{p_i}y_i,T^{q_i}x)+\rho (T^{q_i}x,z)
$$
we get that $\lim_{i\to\infty}\rho (T^{q_i}y,z)=\lim_{i\to\infty}\rho (T^{p_i}y_i,z)=0$.

\ref{viii} By \ref{vi} and the arbitrary choice of $x\in\bigcap_{n=1}^\infty U_n$ it follows that for every $y\in\bigcap_{n=1}^\infty U_n$ there is $z_y\in X$, such that $\lim_{s\to\infty}T^sy=z_y$.
To finish the proof it enough to show that $z_y=z$. Let suppose the contrary i.e.
there exists $y\in\bigcap_{n=1}^\infty U_n$, such that $\lim_{s\to\infty}T^sy=z_y\not= z$.
We have just proved in \ref{vii} that for any $y\in X$ there exists a sequences of natural numbers
$\{q_i\}_{i=1}^\infty$, such that $\{T^{q_i}y\}_{i=1}^\infty$ is convergent to $z$ and consequently $z_y=z$.
\end{proof}

Let for the proofs of the next Corollaries we use the notations from the proof of Theorem \ref{main}.

\begin{Cor}\label{cor1}
Let $(X,\rho )$ be a complete metric space and $T:X\to X$ be a map with the properties:
\renewcommand{\theenumi}{\textnormal{(\alph{enumi}*)}}
\renewcommand{\labelenumi}{\theenumi}
\begin{enumerate}
\item there exist subsets $U_n\subseteq X$, $n\in\N$, such that $T:U_n\to U_n$, and $\bigcap_{n=1}^\infty U_n\not=\emptyset$;
\item there exists $\alpha\in (0,1)$, such that for every $x\in X$ there is $n_0\in\N$, such that for every $u\in \bigcup_{n=n_0}^\infty U_n$ there holds the inequality
$$\rho (Tu,Tx)\leq \alpha\rho (u,x).$$
\end{enumerate}
Then there exists $z\in X$, such that:
\renewcommand{\theenumi}{\textnormal{(\roman{enumi}*)}}
\renewcommand{\labelenumi}{\theenumi}
\begin{enumerate}
\item $Tz=z$;
\item For every $y\in X$ holds $\displaystyle\lim_{i\to\infty}T^{i}y=z$.
\end{enumerate}
\end{Cor}
\begin{proof}
The proof follows immediately by the proof of Theorem \ref{main}, because $n_1(z)=1$, and $p_i=n_1(y_i)=1$.
\end{proof}

\begin{Cor}\label{cor2}
Let $(X,\rho )$ be a complete metric space and $T:X\to X$ be a map with the properties:
\renewcommand{\theenumi}{\textnormal{(\alph{enumi})}}
\renewcommand{\labelenumi}{\theenumi}
\begin{enumerate}
\item there exist subsets $U_n\subseteq X$, $n\in\N$, such that $T:U_n\to U_n$, and $\bigcap_{n=1}^\infty U_n\not=\emptyset$;
\item there exists $\alpha\in (0,1)$, such that for every $x\in X$ there are $n_1, n_0\in\N$, such that for every $u\in \bigcup_{n=n_0}^\infty U_n$ there holds the inequality
$$\rho (T^{n_1}u,T^{n_1}x)\leq \alpha\rho (u,x).$$
\end{enumerate}
Let there holds one of the following conditions:
\renewcommand{\theenumi}{\textnormal{(c.\arabic{enumi})}}
\renewcommand{\labelenumi}{\theenumi}
\begin{enumerate}
\item If there is $x\in\bigcap_{n=1}^\infty U_n$, such that $T$ is continuous at $\lim_{i\to\infty}T^ix$;\label{c1}
\item If there is $x\in\bigcap_{n=1}^\infty U_n$, such that $\lim_{i\to\infty}T^ix\in \bigcap_{n=1}^\infty U_n$;\label{c2}
\item If the sets $\{U_n\}_{n=1}^\infty$ are closed,\label{c3}
\end{enumerate}
then there exists $z\in X$, such that $Tz=z$.
\end{Cor}
\begin{proof}
It follows by Theorem \ref{main} that there exists $z\in X$, such that for any $u\in\bigcap_{n=1}^\infty U_n$ there holds $\lim_{i\to\infty}T^iu=z$ and there exists $n_1(z)$, such that $T^{n_1(z)}=z$.

Let holds \ref{c1}, i.e. $T$ be continuous at $z$. Choose arbitrary $u\in \bigcap_{n=1}^\infty U_n$ and put
$u_i=T^iu$, then for every $\varepsilon >0$ there exists $k_0=k_0(u)\in\N$, such that for every $i\geq k_0$
there hold $\rho (z,Tu_i)=\rho (z,T^iu)<\varepsilon /2$ and $\rho (Tu_i,Tz)<\varepsilon /2$. Now by the inequality
$$
\rho (z,Tz)\leq\rho (z,Tu_i)+\rho (Tu_i,Tz)<\varepsilon
$$
we obtain that $Tz=z$.

Let holds \ref{c2}. Suppose that $Tz=v\not=z$. Then $T^{n_1(z)}v=T^{n_1(z)}(Tz)=T(T^{n_1(z)}z)=Tz=v$.
By the choice of $n_1(z)\in\N$ and Theorem \ref{main} we get
$$
\rho (v,z)=\rho (T^{n_1(z)}v,T^{n_1(z)}z)\leq \alpha\rho (v,z),
$$
which is a contradiction. Consequently $Tz=z$.

Let holds \ref{c3}. The proof follows form \ref{c2}, because if $\{U_n\}_{n=1}^\infty$ are closed then for any
$u\in\bigcap_{n=1}^\infty U_n$ will hold $\lim_{i\to\infty}T^iu\in\bigcap_{n=1}^\infty U_n$.
\end{proof}

\begin{Cor}\label{cor3}
Let $(X,\|\cdot\|)$ be a Banach space and $T:X\to X$ be a map with the properties:
\begin{enumerate}[\quad\textnormal{(\alph{enumi}**)}]
\item there exist subsets $U_n\subseteq X$, $n\in\N$, such that $T:U_n\to U_n$, and $\bigcap_{n=1}^\infty U_n\not=\emptyset$;
\item there exists $\alpha\in (0,1)$, such that for every $x\in X$ there are $n_1, n_0\in\N$, such that for every $u\in \bigcup_{n=n_0}^\infty U_n$ there holds the inequality $$\|T^{n_1}u-T^{n_1}x\|\leq \alpha\|u-x\|.$$
\end{enumerate}
If there is $x\in\bigcap_{n=1}^\infty U_n$, such that the sequence $\{T^ix\}_{i=1}^\infty$ is weakly convergent to
$T\left(\lim_{i\to\infty}T^ix\right)$, then there exists $z\in X$, such that $Tz=z$.
\end{Cor}

\begin{proof}
By Theorem \ref{main} there is $z\in X$, such that for any $x\in\bigcap_{n=1}^\infty U_n$ we have $\lim_{i\to\infty}T^ix=z$.
Let take $x\in\bigcap_{n=1}^\infty U_n$, such that $\{T^ix\}_{i=1}^\infty$ is weakly convergent to $T\left(\lim_{i\to\infty}T^ix\right)=Tz$.
Then $\{T^ix-z\}_{i=1}^\infty$ is weakly convergent to $Tz-z$ and therefore we have
$\|Tz-z\|\leq\underline{\lim}_{i\to\infty}\|T^ix-z\|$.
Now using that $T^ix$ is convergent to $z$ we get that $Tz=z$.
\end{proof}

\section{Examples}

Now we will illustrate the above results with some examples.

{\it Example 1:} Let $f:[0,1]\to [0,1]$ be a convex function, such that $\frac{2}{3}x<f(x)<\frac{4}{5}x$.
We will construct inductively sequences of real numbers $\{a_n\}_{n=1}^\infty$ and $\{\beta_{2n+1}\}_{n=1}^\infty$. Let
$a_1=1$. Put $a_{2n}$, $n\in\N$ to be the solution of the equation
\begin{equation}\label{eq-3}
-a_{2n}+a_{2n-1}=f(a_{2n}),
\end{equation}
$\beta_{2n+1}$, $n\in\N$ to be the solution of the equation
\begin{equation}\label{eq-4}
2a_{2n}+\beta_{2n+1}=f(a_{2n})
\end{equation}
and $a_{2n+1}$, $n\in\N$ to be the solution of the equation
\begin{equation}\label{eq-5}
2a_{2n+1}+f(a_{2n})=2a_{2n}.
\end{equation}
Define $T:[0,1]\to [0,1]$ by
$$
Tx=\left\{
\begin{array}{ll}
-x+a_{2n-1},&n\in\N \quad x\in [a_{2n},a_{2n-1}]\\[12pt]
2x+\beta_{2n+1},&n\in\N\quad x\in [a_{2n+1},a_{2n}]
\end{array}
\right.
$$
and $T(0)=0$. The graphic of the function $T$ is shown in figure 1.
\begin{center}
\epsfig{file=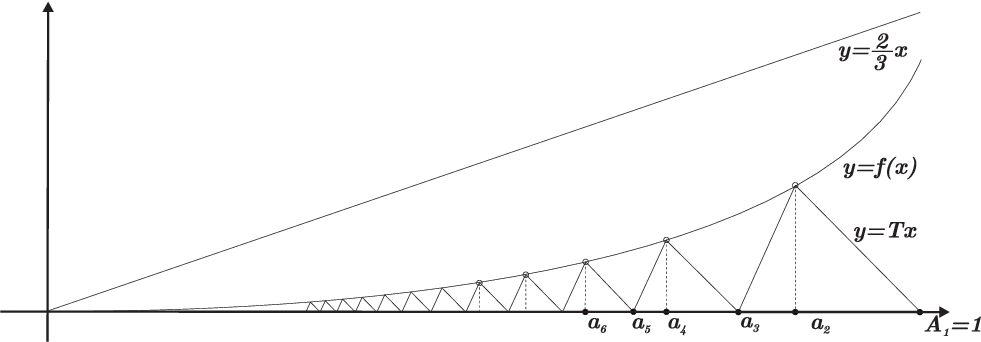,width=12cm}
figure 1
\end{center}

It is easy to see that $\lim_{n\to\infty}a_n=0$. Indeed by \eqref{eq-3} and \eqref{eq-5} we have the inequalities $a_{2n}=a_{2n-1}-f(a_{2n})\leq a_{2n-1}-\frac{2}{3}a_{2n}$ and
$2a_{2n+1}=2a_{2n}-f(a_{2n})\leq 2a_{2n}-\frac{2}{3}a_{2n}=\frac{4}{3}a_{2n}$. After combining the last two inequalities we get that $a_{2n+1}\leq\frac{2}{5}a_{2n-1}$ holds for every $n\in\N$.
By the inequalities $a_{2n+1}<a_{2n}<a_{2n-1}$ it follows that $\lim_{n\to\infty}a_n=0$. Put $U_n=[0,a_{2n+1}]$ for $n\in\N$, then $\bigcap_{n=1}^\infty U_n=\{0\}\not=\emptyset$.

Let $x_0\in [0,1]$ be arbitrary chosen and fixed. Then there is $k_0\in\N$, such that
$x\in [a_{2k_0+1},a_{2k_0-1}]$. By $\lim_{n\to\infty}a_n=0$ it follows that there is $n_0\in\N$, such that
the inequalities $T(u)<T(x_0)$ and $\frac{4}{5}u\leq \frac{4}{5}a_{2n+1}<\frac{4}{5}x_0-T(x_0)$ hold for every $u\in \bigcup_{n=n_0}^\infty U_n$, $n\geq n_1$. Thus for any $u\in \bigcup_{n=n_0}^\infty U_n$ there hold the inequalities
$$
|T(x_0)-T(u)|=T(x_0)-T(u)\leq \frac{4}{5}(x_0-u)-Tu<\frac{4}{5}|x_0-u|
$$
and therefore by Corollary \ref{cor2} it follows that $T$ has a fixed point $z$ and for every $x\in [0,1]$
the sequence $T^nx$ converges to $z$.

It is interesting in this example that there are open sets $V_n=(u_n,t_n)$, such that $\lim_{n\to\infty}u_n=0$ and for any $v_n\in V_n$ there holds the inequality $|T^2v_n-Tv_n|>|Tv_n-v_n|$. Indeed
take $v_n=a_{2n-1}-a_{2n+1}$. By using the inequalities \eqref{eq-3} and \eqref{eq-5} we get
$$
a_{2n}-v_n=2a_{2n}-a_{2n-1}+a_{2n+1}=a_{2n}-\frac{3}{2}f(a_{2n})<0
$$
and therefore $v_n\in (a_{2n},a_{2n-1})$. Then $T(v_n)=-v_n+a_{2n-1}=a_{2n+1}$ and $T^2v_n=Ta_{2n+1}=0$.
Thus $|T^2v_n-Tv_n|=a_{2n+1}$ and $|Tv_n-v_n|=a_{2n-1}-2a_{2n+1}$, because $a_{2n-1}-2a_{2n+1}>0$.
Consequently
$$
\displaystyle\frac{|T^2v_n-Tv_n|}{|Tv_n-v_n|}=\displaystyle\frac{a_{2n+1}}{a_{2n-1}-2a_{2n+1}}=
\displaystyle\frac{a_{2n}-\displaystyle\frac{f(a_{2n})}{2}}{2f(a_{2n})-a_{2n}}>1.
$$
By the continuity of $T$ it follows that there are open sets
$V_n$, such that $v_n\in V_n$ and for any $v\in V_n$ there holds
$$
\displaystyle\frac{|T^2v-Tv|}{|Tv-v|}>1.
$$
{\it Example 2a:} Let us consider the Cantor set.
It is is obtained by repeatedly deleting the open middle thirds of a set of
line segments (starting from $[0,1]$). More specifically, let $K_1$ be the unit interval $[0,1]$
with its middle third   removed, that is $K_1=[0,1/3]\cup [2/3,1]$.
Let $K_2$ be $K_1$ with its middle thirds removed, that is $K_2=[0,1/9]\cup [2/9,3/9]\cup [6/9,7/9]\cup [8/9,1]$.
Continuing in this manner, we generate a sequence of closed sets $K_n$. The Cantor set is defined by $C=\bigcap_{n=1}^\infty K_n\not=\emptyset$.
There is and explicit formula for the open middle third sets that are removed. Put
$V_0=(1/3,2/3)\bigcup (1/9,2/9)\bigcup (7/9,8/9)$
and
$$
V=V_0\bigcup\left(\bigcup_{m=2}^\infty\bigcup_{k=-1}^{m-3}\left(\frac{3^{k+3}-8}{3^m},\frac{3^{k+3}-7}{3^m}\right)\bigcup
\left(\frac{3^{k+3}-2}{3^m},\frac{3^{k+3}-1}{3^m}\right)\right) ,
$$
then $C=[0,1]\backslash V$.
Put $A=V_0\backslash\Q$ and $B=V\backslash \bigl(V_0\cup\Q\bigr)$ and let define
$T_1:A\to\R$ and $T_2:B\to\R$ by
$$
T_1(x)=\left\{
\begin{array}{ll}
\frac{1}{2}x-\frac{1}{18},& x\in \left(\frac{1}{3},\frac{4}{9}\right)\backslash\Q\\[12pt]
\frac{1}{4}x+\frac{1}{18},& x\in \left(\frac{4}{9},\frac{2}{3}\right)\backslash\Q\\[12pt]
\frac{1}{2}x-\frac{1}{54},& x\in \left(\frac{1}{9},\frac{4}{27}\right)\backslash\Q\\[12pt]
\frac{1}{4}x+\frac{1}{54},& x\in \left(\frac{4}{27},\frac{2}{9}\right)\backslash\Q\\[12pt]
\frac{1}{2}x-\frac{7}{54},& x\in \left(\frac{7}{9},\frac{22}{27}\right)\backslash\Q\\[12pt]
\frac{1}{4}x+\frac{2}{27},& x\in \left(\frac{22}{27},\frac{8}{9}\right)\backslash\Q\\[12pt]
\end{array}
\right.
$$
$$
T_2(x)=\left\{
\begin{array}{ll}
\frac{1}{2}x-\frac{3^{k+3}-8}{2.3^{m+1}},& x\in \left(\frac{3^{k+3}-8}{3^m},\frac{3^{k+4}-23}{3^{m+1}}\right)\backslash\Q\\[12pt]
\frac{1}{4}x+\frac{3^{k+3}-7}{4.3^{m+1}},& x\in \left(\frac{3^{k+4}-23}{3^{m+1}},\frac{3^{k+3}-7}{3^m}\right)\backslash\Q\\[12pt]
\frac{1}{2}x-\frac{3^{k+3}-2}{2.3^{m+1}},& x\in \left(\frac{3^{k+3}-2}{3^m},\frac{3^{k+4}-5}{3^{m+1}}\right)\backslash\Q\\[12pt]
\frac{1}{4}x+\frac{3^{k+3}-1}{4.3^{m+1}},& x\in \left(\frac{3^{k+4}-5}{3^{m+1}},\frac{3^{k+3}-1}{3^m}\right)\backslash\Q .\\[12pt]
\end{array}
\right.
$$
Let define a map $T[0,1]\to [0,1]$ by
$$
Tx=\left\{
\begin{array}{ll}
\displaystyle\frac{1}{3}x,&x\in(\Q\cap V)\cup C\\[12pt]
T_1(x),& x\in A\\[12pt]
T_2(x),& x\in B.
\end{array}
\right.
$$
An approximation of the graphic of $T$ is plotted in figure 2.
\begin{center}
\epsfig{file=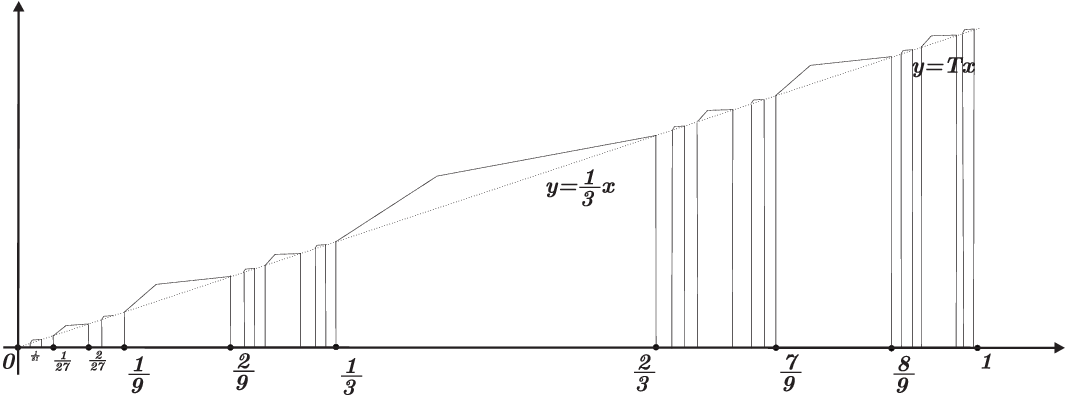,width=12cm}
figure 2
\end{center}
We will show that the map $T$ with the sequence of sets $\{U_n\}_{n=1}^\infty$, where $U_n=K_n$, satisfies the conditions in Theorem \ref{main}.

Let $x_0\in [0,1]$ be arbitrary chosen. There are two cases $x_0\in V$ or $x_0\in C$.\\
Case I) Let $x_0\in V$.

If $x_0\in V\backslash V_0$, then there are $m_0\geq 2$ and $k_0\in\{-1,0,1,\dots m_0-3\}$, such that
$x_0\in \left(\frac{3^{k+3}-8}{3^m},\frac{3^{k+3}-7}{3^m}\right)\bigcup
\left(\frac{3^{k+3}-2}{3^m},\frac{3^{k+3}-1}{3^m}\right)$. Put
$$
r=\sup\left\{\lambda >0:(x_0-\lambda ,x_0+\lambda ) \subset \left(\textstyle\frac{3^{k+3}-8}{3^m},\frac{3^{k+3}-7}{3^m}\right)\cup
\left(\textstyle\frac{3^{k+3}-2}{3^m},\textstyle\frac{3^{k+3}-1}{3^m}\right)\right\}
$$
Choose $n_0=m_0+1$ and $n_1\in\N$ be such that $(1/2)^{n_1-1}<(r/3)$.
Then for every $u\in\bigcup_{n=n_0}^\infty U_n$ there hold the inequalities
$$
|T^{n_1}x_0-T^{n_1}u|\leq \left(\frac{1}{2}\right)^{n_1}x_0+\left(\frac{1}{2}\right)^{n_1}u\leq\left(\frac{1}{2}\right)^{n_1-1}
<\frac{r}{3}<\frac{1}{3}|x_0-u|,
$$
because $Tx\leq\displaystyle\frac{x}{2}$, for every $x\in [0,1]$.

If $x_0\in V_0$, then choose $n_0=2$ and put $r=\sup\left\{\lambda >0:(x_0-\lambda ,x_0+\lambda )\subset V_0\right\}$.
Choose $n_1$ be such that $(1/2)^{n_1-1}<(r/3)$.
Then for every $u\in\bigcup_{n=2}^\infty U_n$ there hold the inequalities
$$
|T^{n_1}x_0-T^{n_1}u|\leq \left(\frac{1}{2}\right)^{n_1}x_0+\left(\frac{1}{2}\right)^{n_1}u\leq\left(\frac{1}{2}\right)^{n_1-1}
<\frac{r}{3}<\frac{1}{3}|x_0-u|.
$$
Case II)
Let $x_0\in C$. Then for any $u\in \bigcup_{n=2}^\infty U_n$ there holds
$$
|Tx_0-Tu|\leq \frac{1}{3}|x_0-u|,
$$
because $Tx\geq\displaystyle\frac{x}{3}$ for every $x\in [0,1]$ and $Tx=\displaystyle\frac{x}{3}$ for every $x\in C$. Now we can apply Theorem \ref{main}.
Let us mention the result
\begin{theo}\label{th-S}\cite{S}
Let $X$ be a Banach space, and $T:X\to X$ a continuous mapping satisfying the condition: there exists a constant
$\alpha\in (0,1)$ such that for each $x\in X$, there is a positive integer $n(x)$ such that for all $y\in X$
$$
\rho (T^{n(x)}y,T^{n(x)}x)\leq \alpha\rho (y,x).
$$
Then $T$ has a unique fixed point $z$ and $\displaystyle\lim_{s\to\infty}T^sx=z$ for each $x\in X$.
\end{theo}
It is easy to see that for any $x\in V\cap\Q$ and for any $n_1\in\N$ we can choose an irrational $u$, such that
$|T^{n_1}x-T^{n_1}u|>|x-u|$ and thus the conditions in Theorem \ref{th-S}
are not satisfied.\\
{\it Example 2b:} Let $T:[0,1]\to [0,1]$ be defined as in Example 2a, for every $x\in (0,1)\backslash \{1/2\}$. Let $T(0)=1$, $T(1)=1/2$, $T(1/2)=1$. Then all he conditions in Theorem \ref{main} are satisfied with with the sets $U_n=[0,1/n]$, $n\in\N$. It easy to see that $n_1(0)=3$. Therefore $T^3(0)=0$ and obviously there is no $x\in [0,1]$, such that $Tx=x$.

{\it Example 3:} Let $(X,\|\cdot\|)$ be a Banach space with a basis $\{e_i\}_{i=1}^\infty$,
such that
\begin{equation}\label{eq-6}
\begin{array}{l}
\mbox{if}\quad x=\displaystyle\sum_{i=1}^\infty x_ie_i, y=\displaystyle\sum_{i=1}^\infty y_ie_i\quad \mbox{with}\quad
|x_i|\leq |y_i|\quad \mbox{for every}\quad i\in\N,\\
\mbox{then}\quad \|x\|\leq\|y\|.
\end{array}
\end{equation}

Let $\{\alpha_i\}_{i=1}^\infty$ be an increasing sequence of positive reals, convergent to $1$. Let
$T:X\to X$, be linear map defined by $Te_k=\alpha_ke_k$.
Consider the sets $U_n=\left\{\sum_{i=1}^n\lambda_ie_i:\left\|\sum_{i=1}^n\lambda_ie_i\right\|\leq \frac{1}{n}\right\}$.
Obviously $\bigcap_{n=1}^\infty U_n\not=\emptyset$.

Let $x=\sum_{i=1}^\infty x_ie_i\in X$ be arbitrary chosen. There is $m_1\in\N$, such that
$$
\left\|\displaystyle\sum_{i=m_1}^\infty x_ie_i\right\|<\frac{1}{8}\|x\|,
$$
which implies the inequality $\frac{7}{8}\|x\|\leq\left\|\sum_{i=1}^{m_1}x_ie_i\right\|$.
Choose $m_2\in\N$, such that $\frac{1}{m_2}\leq\frac{1}{8}\|x\|$ and put $m_0=\max\{m_1,m_2\}$.
Then for every $m\geq m_0$ and every $y\in U_m$ there hold the inequalities
$$
\begin{array}{lll}
\left\|\displaystyle\sum_{i=m+1}^\infty x_ie_i\right\|&\leq&
\displaystyle\frac{1}{8}\|x\|\leq\displaystyle\frac{1}{4}\left(\displaystyle\frac{3}{4}\|x\|-\displaystyle\frac{1}{m}\right)\leq\displaystyle\frac{1}{4}\left|\left\|\displaystyle\sum_{i=1}^m x_ie_i\right\|-\|y\|\right|\\[12pt]
&\leq&\displaystyle\frac{1}{4}\left|\|x\|-\|y\|\right|\leq\displaystyle\frac{1}{4}\|x-y\|
\end{array}
$$
and
$$
\|y\|\leq \frac{1}{m}\leq\frac{1}{8}\|x\|\leq\frac{1}{4}\|x\|-\frac{1}{4}\|y\|\leq\frac{1}{4}\|x-y\|.
$$
Choose $s_0=s(x)$, so that $\max\{\alpha_i^{s_0}\leq 2^{-2}:i=1,2,\dots ,m_0\}$. Let $y\in\bigcup_{n=m_0}^\infty U_n$, then there is $n\geq m_0$, such that $y\in U_n$ and $y=\sum_{i=1}^\infty\lambda_ie_1$, where $\lambda_i=0$ for every $i>n$. We can write the chain of inequalities:
$$
\begin{array}{lll}
\left\|T^{s_0}x-T^{s_0}y\right\|&=&\left\|\displaystyle\sum_{i=1}^{m_0}\alpha_i^{s_0}(x_i-\lambda_i)e_i+\displaystyle\sum_{i=m_0+1}^{\infty}\alpha_i^{s_0}x_ie_i\right\|\\[12pt]
&\leq&\left\|\displaystyle\sum_{i=1}^{m_0}\alpha_i^{s_0}(x_i-\lambda_i)e_i+\displaystyle\sum_{i=m_0+1}^{\infty}\alpha_{m_0}^{s_0}(x_i-\lambda_i)e_i\right\|\\[12pt]
&&+\left\|\displaystyle\sum_{i=m_0+1}^{\infty}(\alpha_i^{s_0}-\alpha_{m_0}^{s_0})x_ie_i\right\|+\left\|\displaystyle\sum_{i=m_0+1}^{\infty}(\alpha_i^{s_0}-\alpha_{m_0}^{s_0})\lambda_ie_i\right\|\\[12pt]
&\leq&\displaystyle\frac{1}{4}\|x-y\|+\left\|\displaystyle\sum_{i=m_0+1}^{\infty}x_ie_i\right\|+\|y\|\\[12pt]
&\leq&\displaystyle\frac{1}{4}\|x-y\|+\frac{1}{4}\|x-y\|+\frac{1}{4}\|x-y\|=\displaystyle\frac{3}{4}\|x-y\|.
\end{array}
$$

The sets $U_n$ are closed for every $n\in\N$ and by Corollary \ref{cor2} there exists $z\in X$, such that $Tz=z$.  It is easy to see that $\lim_{s\to\infty}T^se_1=\sigma$ and $\sigma$ is the fixed point of the map $T$, where $\sigma$ is the zero vector in $X$. By Theorem \ref{main} we have that for every $y\in X$ there is a sequence of naturals $\{q_i\}_{i=1}^\infty$, such that $\lim_{i\to\infty}T^{q_i}y=\sigma$.
Therefore for every $\varepsilon >0$, there exists $k_\varepsilon\in\N$, such that $\|T^{q_{k_\varepsilon}}y\|<\varepsilon$.
By the fact that $0<\alpha_i<1$ for every $i\in\N$ it follows that for every $s\geq q_{k_\varepsilon}$ holds
$|\alpha_i|^s\leq |\alpha_i|^{q_{k_\varepsilon}}$ for every $i\in\N$. Now by \eqref{eq-6} it follows that for every $s\geq q_{k_\varepsilon}$ there holds the inequality $\|T^sy\|<\varepsilon$. By the arbitrary choice of $\varepsilon >0$ we get that $\displaystyle\lim_{s\to\infty}T^sy=\sigma$.

A Result that generalizes Theorem \ref{th-S} is the following:

\begin{theo}\label{th-G}\textnormal{\cite{G}}
Let $X$ be a Banach space, and let $T:X\to X$ a mapping. Suppose there exists $B\subset X$ such that
\begin{enumerate}[\textnormal{(G\arabic{enumi})}]
  \item $T(B)\subset B$;
  \item for some $\alpha\in (0,1)$ and each $y\in B$ there is an integer $n(y)\geq 1$ with
$$
\rho (T^{n(y)}x,T^{n(y)}y)\leq \alpha\rho (x,y)
$$
for all $x\in B$;
  \item for some $x\in B$, ${\mathrm{ cl}}\{T^sx:s\in\N\}\subset B$.
\end{enumerate}
Then there is a unique $z\in B$ such that $T(z)=z$ and $\displaystyle\lim_{s\to\infty}T^sx=z$ for each $x\in B$. Furthermore, if
\begin{equation}\label{eq-7}
\rho (T^{n(z)}x,T^{n(z)}z)\leq \alpha\rho (x,z)
\end{equation}
for all $x\in X$, then $z$ is unique in $X$ and $\displaystyle\lim_{s\to\infty}T^sx=z$ for each $x\in X$.
\end{theo}

It is easy to see in Example 3, that for every $0<\alpha <1$ and every $n\in\N$ there is $s_n\in\N$
such that
$$
\|T^ne_{s_n}-T^n\sigma\|=\|T^ne_{s_n}\|=(\alpha_{s_n})^n\|e_{s_n}\|\geq \alpha\|e_{s_n}\|=\alpha\|e_{s_n}-\sigma\|,
$$
which shows that the condition \eqref{eq-7} in Theorem \ref{th-G} is not satisfied.

\end{document}